\documentclass[12pt]{article}%
\usepackage{graphicx}
\usepackage[intlimits]{amsmath}
\usepackage{latexsym}
\usepackage{amsfonts}
\usepackage{amssymb}%
\makeatletter
\newfont{\footsc}{cmcsc10 at 8truept}
\newfont{\footbf}{cmbx10 at 8truept}
\newfont{\footrm}{cmr10 at 10truept}
\newtheorem{theorem}{Theorem}

\newtheorem{corollary}[theorem]{Corollary}

\newtheorem{lemma}[theorem]{Lemma}

\newenvironment{proof}[1][Proof]{\noindent{\textbf {#1}  }}  {\hfill$\Box$\bigskip}

\ifx\pdfoutput\relax\let\pdfoutput=\undefined\fi
\newcount\msipdfoutput
\ifx\pdfoutput\undefined\else
\ifcase\pdfoutput\else
\ifx\paperwidth\undefined\else
\ifdim\paperheight=0pt\relax\else\pdfpageheight\paperheight\fi
\ifdim\paperwidth=0pt\relax\else\pdfpagewidth\paperwidth\fi
\fi\fi\fi
\begin{document}

\title{Cycles and Stability\\(full version)}
\author{V. Nikiforov and R. H. Schelp\\Department of Mathematical Sciences, University of Memphis, \\Memphis, Tennessee, 38152\\{\small E-mail address: V. Nikiforov, \textit{vnikifrv@memphis.edu;}\ R.H.
Schelp, \textit{rschelp@memphis.edu}}}
\maketitle

\begin{abstract}
We prove a number of Tur\'{a}n and Ramsey type stability results for cycles,
in particular, the following one:

Let $n\geq4$, $0<\beta\leq1/2-1/2n$, and the edges of $K_{\left\lfloor \left(
2-\beta\right)  n\right\rfloor }$ be $2$-colored so that no mononchromatic
$C_{n}$ exists. Then, for some $q\in\left(  \left(  1-\beta\right)
n-1,n\right)  ,$ we may drop a vertex $v$ so that in $K_{\left\lfloor \left(
2-\beta\right)  n\right\rfloor }-v$ one of the colors induces
$K_{q,\left\lfloor \left(  2-\beta\right)  n\right\rfloor -q-1},$ while the
other one induces $K_{q}\cup K_{\left\lfloor \left(  2-\beta\right)
n\right\rfloor -q-1}.$

We also derive the following Ramsey type result.

If $n$ is sufficiently large and $G$ is a graph of order $2n-1,$ with minimum
degree $\delta\left(  G\right)  \geq\left(  2-10^{-6}\right)  n,$ then for
every $2$-coloring of $E\left(  G\right)  $ one of the colors contains cycles
$C_{t}$ for all $t\in\left[  3,n\right]  .$

\end{abstract}

\section{Introduction}

Our graph theoretic notation is standard (e.g., see \cite{Bol98}). We write
$c\left(  G\right)  $ for the length of the longest cycle in a graph$\ G.$
Given a graph $G$ and disjoint sets $A,B\subset V\left(  G\right)  $, we write
$E_{G}\left(  A,B\right)  $ for the set of $A-B$ edges of $G,$ and, abusing
notation, $A\times B$ for the set of all possible $A-B$ edges. $A^{\left(
k\right)  }$ stands for the family of $k$-subsets of a set $A.$

Given a graph $G,$ a\emph{ }$2$\emph{-coloring} of $E\left(  G\right)  $ is a
partition $E\left(  G\right)  =E\left(  R\right)  \cup E\left(  B\right)  ,$
where $R$ and $B$ are graphs with $V\left(  R\right)  =V\left(  B\right)
=V\left(  G\right)  .$ Given a $2$-coloring $E\left(  G\right)  =E\left(
R\right)  \cup E\left(  B\right)  ,$ a statement $S\left(  R,B\right)  $
involving $R$ and $B$ is said to be true \emph{up to color,} if either
$S\left(  R,B\right)  $ or $S\left(  B,R\right)  $ is true.

In this paper we study stability results about cycles. Stability is a topic
studied mostly in extremal problems of Tur\'{a}n type, but also appears neatly
in Ramsey problems for cycles, as shown below, and for books, as shown in
\cite{NRS05}. In addition to being interesting for their own sake, stability
theorems are excellent tools for obtaining new extremal results, like Theorem
\ref{arrth} below. In fact, our initial motivation came from another
application that will be presented in a sequel to this paper. It should be
noted that our most applicable results - Theorems \ref{thDC} and \ref{th3par}
- are too technical to be stated in the Introduction.

Our first stability result is of Tur\'{a}n type.

\begin{theorem}
\label{cycth}Let $0<\gamma<10^{-5}$. If $G=G\left(  n\right)  $ is a graph
with $e\left(  G\right)  >n^{2}/4,$ then one of the following two conditions hold:

(a) $C_{t}\subset G$ for every $t\in\left[  3,\left\lceil \left(
1/2+\gamma\right)  n\right\rceil \right]  ;$

(b) there exists a vertex $v$ such that $G-v=H_{1}\cup H_{2},$ where $H_{1}$
and $H_{2}$ are vertex-disjoint graphs satisfying
\[
\left(  \frac{1}{2}-900\gamma\right)  n<\left\vert H_{1}\right\vert
\leq\left\vert H_{2}\right\vert <\left(  \frac{1}{2}+900\gamma\right)  n.
\]

\end{theorem}

The following theorem is a stability result of Ramsey type. It states that if
$p$ is close to $2n,$ and the edges of $K_{p}$ are $2$-colored so that no
monochromatic cycle $C_{n}$ exists, then we may remove a vertex from $K_{p}$
so that, for some $q$ close to $n,$ one of the colors induces $K_{q,p-q-1},$
while the other one induces $K_{q}\cup K_{p-q-1}$.

\begin{theorem}
\label{cycth1} Let $n\geq4,$ $0<\beta\leq\left\lfloor n/2\right\rfloor /n$,
and $E\left(  K_{\left\lfloor \left(  2-\beta\right)  n\right\rfloor }\right)
=E\left(  R\right)  \cup E\left(  B\right)  $ be a $2$-coloring such that
$C_{n}\nsubseteq R$ and $C_{n}\nsubseteq B.$ Then there exist a vertex $u\in
V\left(  K_{\left\lfloor \left(  2-\beta\right)  n\right\rfloor }\right)  $
and a partition $V\left(  K_{\left\lfloor \left(  2-\beta\right)
n\right\rfloor }\right)  =U_{1}\cup U_{2}\cup\left\{  u\right\}  $ with
\begin{equation}
\left(  1-\beta\right)  n-1<\left\vert U_{1}\right\vert \leq\left\vert
U_{2}\right\vert <n \label{thin}%
\end{equation}
satisfying, up to color,
\[
E\left(  R-v\right)  =U_{1}^{\left(  2\right)  }\cup U_{2}^{\left(  2\right)
}\text{ \ \ and \ \ }E\left(  B-v\right)  =U_{1}\times U_{2}.
\]

\end{theorem}

We also derive the following Ramsey type result.

\begin{theorem}
\label{arrth}If $n$ is sufficiently large and $G=G\left(  2n-1\right)  $ is a
graph with $\delta\left(  G\right)  \geq\left(  2-10^{-6}\right)  n,$ then,
for every $2$-coloring $E\left(  G\right)  =E\left(  R\right)  \cup E\left(
B\right)  ,$ either $C_{t}\subset R$ for all $t\in\left[  3,n\right]  $ or
$C_{t}\subset B$ for all $t\in\left[  3,n\right]  .$
\end{theorem}

The rest of the paper is organized as follows: in Section \ref{Suc} we give
sufficient conditions for cycles and paths, in Section \ref{Ext} we prove
Tur\'{a}n type stability results including Theorem \ref{cycth}, and in Section
\ref{Ram} we prove Ramsey type stability results including Theorem
\ref{cycth1} and Theorem \ref{arrth}.

\section{\label{Suc}Sufficient conditions for cycles and paths}

In this section we list sufficient conditions for the existence of cycles and
paths. Most of the them are known, but we also give a few new ones.

\begin{theorem}
[Erd\H{o}s and Gallai \cite{ErGa59}]\label{leEG}If $e\left(  G\right)
\geq\left\vert G\right\vert $, then $c\left(  G\right)  >2e\left(  G\right)
/\left\vert G\right\vert .$\hfill$\square$
\end{theorem}

This result was significantly improved for $2$-connected graphs by Woodall
\cite{Woo76}, and recently by Fan, Lv, and Weng \cite{FLW04}.

\begin{theorem}
[Fan, Lv, and Weng \cite{FLW04}]\label{ThWo}If the graph $G=G\left(  n\right)
$ is $2$-connected and $c\left(  G\right)  =c$ then%
\[
e\left(  G\right)  \leq\binom{c+1-\left\lfloor c/2\right\rfloor }%
{2}+\left\lfloor \frac{c}{2}\right\rfloor \left(  n-c-1+\left\lfloor \frac
{c}{2}\right\rfloor \right)  .
\]
\hfill$\square$
\end{theorem}

\begin{theorem}
[Bollob\'{a}s \cite{Bol78}, p. 150 ]\label{ThBe}If $G$ is a graph with
$e\left(  G\right)  >\left\vert G\right\vert ^{2}/4$, then $C_{t}\subset G$
for every $3\leq t\leq c\left(  G\right)  .$\hfill$\square$
\end{theorem}

Implicit in \cite{ErGa59} (see \cite{Bon95}, p. 26) is the following theorem.

\begin{theorem}
[Erd\H{o}s and Gallai \cite{ErGa59}]\label{ThEG} If $G$ is a graph with
$\delta\left(  G\right)  >\left\vert G\right\vert /2$, then every two vertices
of $G$ can be joined by a path of order $\left\vert G\right\vert $%
.\hfill$\square$
\end{theorem}

The following theorem follows from results of Brandt, Faudree, and Goddard
(\cite{BFG98}, p. 143).

\begin{theorem}
[Brandt, Faudree, and Goddard \cite{BFG98}]\label{thBFG}If $n>30$ and
$G=G\left(  n\right)  $ is a $2$-connected, nonbipartite graph with
$\delta\left(  G\right)  >2n/5,$ then $C_{t}\subset G$ for all $t\in\left[
3,c\left(  G\right)  \right]  .$\hfill$\square$
\end{theorem}

We derive below a simple consequence of Theorem \ref{ThWo}.

\begin{corollary}
\label{CorWo}If $G=G\left(  n\right)  $ is a $2$-connected graph, then either
$G$ is Hamiltonian or%
\[
c\left(  G\right)  >2n\left(  1-\sqrt{1-\frac{2e\left(  G\right)  }{n^{2}}%
}\right)  .
\]

\end{corollary}

\begin{proof}
Write $m$ for $e\left(  G\right)  $ and let $c\left(  G\right)  $ be even, say
$c\left(  G\right)  =2k.$ Theorem \ref{ThWo} implies that%
\[
m\leq\frac{k\left(  k+1\right)  }{2}+k\left(  n-k-1\right)  =-\frac{k\left(
k+1\right)  }{2}+kn.
\]
Hence $k^{2}-2kn+2m\leq-k<0,$ and the assertion follows.

Let now $c\left(  G\right)  $ be odd, say $c\left(  G\right)  =2k+1.$ Theorem
\ref{ThWo} implies that
\begin{align*}
m  &  \leq\binom{2k+2-k}{2}+k\left(  n-2k-2+k\right)  =\frac{\left(
k+2\right)  \left(  k+1\right)  }{2}+k\left(  n-k-2\right) \\
&  =-\frac{k\left(  k+1\right)  }{2}+kn+1;
\end{align*}
hence, $k^{2}-k\left(  2n-1\right)  +2m-2\leq0.$ In view of $\left(
2n-1\right)  ^{2}-8\left(  m-1\right)  \leq4n^{2}-8m,$ it follows that
\[
2k+1\geq2n-\sqrt{\left(  2n-1\right)  ^{2}-8\left(  m-1\right)  }\geq
2n-\sqrt{4n^{2}-8m},
\]
completing the proof.\bigskip
\end{proof}

We shall make use of the following simple consequence of Theorem \ref{ThEG}.

\begin{lemma}
\label{2hl} If $G=G\left(  n\right)  $ is a graph with $\delta\left(
G\right)  \geq n/2+1,$ every two vertices of $G$ can be joined by paths of
order $n$ and $n-1.$
\end{lemma}

\begin{proof}
Select $u,v\in V\left(  G\right)  .$ Theorem \ref{ThEG} implies that $u$ and
$v$ may be joined by a path of order $n.$ From $n-1\geq\delta\left(  G\right)
\geq n/2+1$ we see that $n\geq4.$ Select $w\neq u,v$ and consider $G^{\prime
}=G-w.$ We have
\[
\delta\left(  G^{\prime}\right)  \geq\delta\left(  G\right)  -1\geq\frac{n}%
{2}>\frac{n-1}{2}=\frac{\left\vert G^{\prime}\right\vert }{2},
\]
thus, again by Theorem \ref{ThEG}, $u$ and $v$ can be joined by a path of
order $n-1.$
\end{proof}

\begin{lemma}
\label{ecl}Let $G$ be a bipartite graph with vertex classes $A$ and $B,$
$\left\vert A\right\vert \leq\left\vert B\right\vert ,$ and $\delta
=\delta\left(  G\right)  \geq\left\vert B\right\vert /2+1.$ Then

(i) if $x,y\in A$ or $x,y\in B,$ $G$ contains an $xy$-path of length $t$ for
all even $t\in\left[  2,2\left(  2\delta-\left\vert A\right\vert -1\right)
\right]  ;$

(ii) if $x\in A$, $y\in B,,$ $G$ contains an $xy$-path of length $t$ for all
odd $t\in\left[  3,2\left(  2\delta-\left\vert A\right\vert -1\right)
\right]  ;$

(iii) $C_{t}\subset G$ for all even $t\in\left[  4,2\left(  2\delta-\left\vert
A\right\vert -1\right)  \right]  .$
\end{lemma}

\begin{proof}
To prove \emph{(i)} and \emph{(ii)} we use induction on $t.$ If $x,y\in A$ or
$x,y\in B$, then
\[
\left\vert \Gamma\left(  x\right)  \cap\Gamma\left(  y\right)  \right\vert
\geq d\left(  x\right)  +d\left(  y\right)  -\left\vert B\right\vert
>2\delta-\left\vert B\right\vert \geq2,
\]
and so, there exists an $xy$-path of length $2$. If $x\in A$, $y\in B,$ select
$u_{1}\in\Gamma\left(  x\right)  .$ Since $\left\vert \Gamma\left(
u_{1}\right)  \cap\Gamma\left(  y\right)  \right\vert \geq2,$ there exists
$u_{2}\in\left(  \Gamma\left(  u_{1}\right)  \cap\Gamma\left(  y\right)
\right)  \backslash\left\{  x\right\}  ;$ the path $x,u_{1},u_{2},y$ has
length $3$. To complete the induction we show that if $l<2\left(
2\delta-\left\vert A\right\vert -1\right)  ,$ every $xy$-path $P=xu_{1}%
,...u_{l-1},y$ of length $l$ can be extended to a $xy$-path of length $l+2.$
Select $u_{i},u_{i+1}\in P$ so that $u_{i}\in A,$ $u_{i+1}\in B.$ Since
\[
\left\vert P\cap B\right\vert \leq\frac{l+1}{2}<2\delta-\left\vert
A\right\vert <\delta,
\]
we can select a vertex $v\in\Gamma\left(  u_{i}\right)  \backslash P.$ Since
\[
\left\vert \Gamma\left(  u_{i+1}\right)  \cap\Gamma\left(  v\right)
\right\vert \geq2\delta-\left\vert A\right\vert >\frac{l+1}{2}\geq\left\vert
P\cap A\right\vert ,
\]
we can select $w\in\Gamma\left(  u_{i+1}\right)  \cap\Gamma\left(  v\right)
\backslash P.$ The $xy$-path
\[
x,u_{1},...u_{i},v,w,u_{i+1},...,y
\]
has length $l+2$, completing the induction and the proof of \emph{(i)} and
\emph{(ii)}.

To prove \emph{(iii)},\emph{ }select two adjacent vertices $x\in A$, $y\in B.$
According to \emph{(ii) }there exists an $xy$-path of odd length $t\in\left[
3,2\left(  2\delta-\left\vert A\right\vert -1\right)  -1\right]  ,$ and
consequently, a cycle of length $t+1,$ completing the proof.
\end{proof}

\section{\label{Ext}Tur\'{a}n type stability}

Most results in this paper are derived from the following theorem.

\begin{theorem}
\label{thDC}Let $0<\alpha<10^{-5},$ $0\leq\beta<10^{-5},$ and $n\geq
\alpha^{-1}/2.$ If $G=G\left(  n\right)  $ is a graph with $e\left(  G\right)
>\left(  1/4-\beta\right)  n^{2}$, then one of the following conditions holds:

(i) $c\left(  G\right)  \geq\left(  1/2+\alpha\right)  n;$

(ii) there exists a set $M\subset V\left(  G\right)  $ with $\left\vert
M\right\vert <840\left(  \alpha+2\beta\right)  n$ such that $G-M$ consists of
two components $G_{1},G_{2}$ satisfying%
\begin{align}
\left(  \frac{1}{2}-840\left(  \alpha+2\beta\right)  \beta\right)  n  &
<\left\vert G_{1}\right\vert \leq\left\vert G_{2}\right\vert <\left(  \frac
{1}{2}+20\left(  \alpha+2\beta\right)  \right)  n,\label{in1}\\
\delta\left(  G_{1}\right)   &  \geq\frac{3n}{7},\text{ \ \ \ \ }\delta\left(
G_{2}\right)  \geq\frac{3n}{7}. \label{in2}%
\end{align}

\end{theorem}

\begin{proof}
Assume that \emph{(i) }fails, i.e.,%
\begin{equation}
c\left(  G\right)  <\left(  1/2+\alpha\right)  n. \label{assum}%
\end{equation}

The rest of our proof has two phases - in the first one we find $M_{1}\subset
V\left(  G\right)  $ such that $\left\vert M_{1}\right\vert <40\left(
\alpha+\beta\right)  n$ and $G-M_{1}$ consists of two components $H_{1},H_{2}$
satisfying%
\begin{equation}
\left(  \frac{1}{2}-20\alpha+40\beta\right)  n<\left\vert H_{1}\right\vert
\leq\left\vert H_{2}\right\vert <\left(  \frac{1}{2}+20\alpha+40\beta\right)
n. \label{in0}%
\end{equation}

Then, in the second phase, we obtain $G_{1}$ and $G_{2}$ by dropping the low
degree vertices from $H_{1}$ and $H_{2}.$

Setting%
\[
M_{0}=\left\{  v:v\in V\left(  G\right)  ,\text{ }d\left(  v\right)
\leq9n/40\right\}  ,
\]
our first goal is to prove that
\begin{equation}
\left\vert M_{0}\right\vert <\left(  20\alpha+40\beta\right)  n. \label{ubm}%
\end{equation}
Indeed, Lemma \ref{leEG} implies that
\[
2e\left(  G-M_{0}\right)  <c\left(  G-M_{0}\right)  \left(  n-\left\vert
M_{0}\right\vert \right)  \leq c\left(  G\right)  \left(  n-\left\vert
M_{0}\right\vert \right)  <\left(  \frac{1}{2}+\alpha\right)  \left(
n-\left\vert M_{0}\right\vert \right)  n,
\]
and so,%
\begin{align*}
\left(  \frac{1}{4}+\frac{\alpha}{2}\right)  n^{2}-\frac{1}{4}\left\vert
M_{0}\right\vert n  &  \geq\frac{1}{2}\left(  \frac{1}{2}+\alpha\right)
\left(  n-\left\vert M_{0}\right\vert \right)  n>e\left(  G-M_{0}\right)  \geq
e\left(  G\right)  -\sum_{u\in M_{0}}d\left(  u\right) \\
&  >\left(  \frac{1}{4}-\beta\right)  n^{2}-\frac{9n}{40}\left\vert
M_{0}\right\vert ,
\end{align*}
implying (\ref{ubm}).

From
\[
\left(  \frac{11}{10}-8\beta\right)  n\geq\left(  20\alpha+40\beta\right)
n\geq\left(  1-4\beta\right)  20\left(  \alpha+2\beta\right)  n\geq\left(
1-4\beta\right)  \left\vert M_{0}\right\vert
\]
we deduce that
\[
\left(  \frac{1}{4}-\beta\right)  n^{2}-\frac{9}{40}\left\vert M_{0}%
\right\vert n\geq\left(  \frac{1}{4}-\beta\right)  \left(  n-\left\vert
M_{0}\right\vert \right)  ^{2}.
\]
If $\kappa\left(  G-M_{0}\right)  \geq2,$ Corollary \ref{CorWo} implies that
\begin{align*}
c\left(  G\right)   &  \geq2\left(  n-\left\vert M_{0}\right\vert \right)
\left(  1-\sqrt{1-\frac{2e\left(  G-M_{0}\right)  }{\left(  n-\left\vert
M_{0}\right\vert \right)  ^{2}}}\right) \\
&  \geq2\left(  1-20\alpha-40\beta\right)  \left(  1-\sqrt{\frac{1}{2}+2\beta
}\right)  n\\
&  \geq2\left(  1-20\alpha-40\beta\right)  \left(  1-\frac{1+2\beta}{\sqrt{2}%
}\right)  n\\
&  \geq\left(  2-\sqrt{2}\right)  \left(  1-20\alpha-40\beta\right)  \left(
1-2\left(  \sqrt{2}+1\right)  \beta\right)  n\\
&  \geq\left(  2-\sqrt{2}\right)  \left(  1-20\alpha-45\beta\right)
n\geq\left(  \frac{1}{2}+\alpha\right)  n,
\end{align*}
contradicting (\ref{assum}).

Hence, there exists $K\subset V\left(  G\right)  $ with $\left\vert
K\right\vert \leq1$ such that the graph $G^{\prime}=G-M_{0}-K$ is
disconnected. Observe that $\alpha+2\beta<3\times10^{-5}$ and $n>10^{5}$ imply%
\begin{align*}
\delta\left(  G^{\prime}\right)   &  =\delta\left(  G-M_{0}-K\right)
>\frac{9n}{40}-\left\vert M_{0}\right\vert -1>\frac{9n}{40}-20\left(
\alpha+2\beta\right)  n-1\\
&  \geq\left(  \frac{9}{40}-\frac{20\times3}{10^{5}}\right)  n-1\geq\frac
{n}{5}.
\end{align*}

\emph{Case 1: }$G^{\prime}$\emph{ has a component }$G^{\prime\prime}$\emph{
with }$\left\vert G^{\prime\prime}\right\vert \leq n/3$

Then, by Lemma \ref{leEG},
\[
2e\left(  G^{\prime}-G^{\prime\prime}\right)  \leq c\left(  G^{\prime
}-G^{\prime\prime}\right)  \left(  \left\vert G^{\prime}\right\vert
-\left\vert G^{\prime\prime}\right\vert \right)  \leq c\left(  G\right)
\left(  n-\left\vert M_{0}\right\vert -\left\vert G^{\prime\prime}\right\vert
\right)  .
\]
In view of $\Delta\left(  G^{\prime\prime}\right)  <n/3,$ we see that
\begin{align*}
\left(  n-\left\vert M_{0}\right\vert -\left\vert G^{\prime\prime}\right\vert
\right)  \left(  \frac{1}{4}+\frac{\alpha}{2}\right)  n  &  \geq e\left(
G^{\prime}-G^{\prime\prime}\right)  \geq e\left(  G^{\prime}\right)  -e\left(
G^{\prime\prime}\right) \\
&  \geq\left(  \frac{1}{4}-\beta\right)  n^{2}-\frac{9n}{40}\left\vert
M_{0}\right\vert -n\left\vert K\right\vert -\frac{\left\vert G^{\prime\prime
}\right\vert n}{6},
\end{align*}
and therefore,%
\[
\left(  \frac{1}{4}+\frac{\alpha}{2}\right)  n-\left(  \frac{1}{4}%
+\frac{\alpha}{2}\right)  \left\vert M_{0}\right\vert -\left(  \frac{1}%
{4}+\frac{\alpha}{2}\right)  \left\vert G^{\prime\prime}\right\vert
\geq\left(  \frac{1}{4}-\beta\right)  n-\frac{9}{40}\left\vert M_{0}%
\right\vert -\left\vert K\right\vert -\frac{\left\vert G^{\prime\prime
}\right\vert }{6},
\]
implying that%
\begin{align*}
\left(  \frac{\alpha+2\beta}{2}\right)  n-\left(  \frac{1}{40}+\frac{\alpha
}{2}\right)  \left\vert M_{0}\right\vert +\left\vert K\right\vert  &
\geq\left(  \frac{1}{12}+\frac{\alpha}{2}\right)  \left\vert G^{\prime\prime
}\right\vert \\
&  >\left(  \frac{1}{12}+\frac{\alpha}{2}\right)  \delta\left(  G^{\prime
}\right)  >\left(  \frac{1}{12}+\frac{\alpha}{2}\right)  \frac{1}{5}n.
\end{align*}
This gives
\[
6\left(  \alpha+2\beta\right)  n+12\geq\frac{1+6\alpha}{5}n,
\]
and, in view of $\alpha<10^{-5},\beta\leq10^{-5},$ it follows that%
\[
n\leq24\alpha n+60\beta n+60<\frac{84}{10^{5}}n+60.
\]
This inequality is a contradiction, as $n\geq\alpha^{-1}/2>10^{5}/2.$

\emph{Case 2: The order of each component of }$G^{\prime}$\emph{ is greater
than }$n/3$

Therefore, $G^{\prime}$ has exactly two components - $H_{1}$ and $H_{2};$ let,
say $\left\vert H_{1}\right\vert \leq\left\vert H_{2}\right\vert .$ Setting
$M_{1}=M_{0}\cup K,$ we see that
\[
\left\vert M_{1}\right\vert \leq20\left(  \alpha+2\beta\right)  n+1\leq
20\left(  \alpha+2\beta\right)  n+2\alpha n<40\left(  \alpha+\beta\right)  n.
\]

We shall prove that inequalities (\ref{in0}) hold. From Lemma \ref{leEG} we
have%
\[
e\left(  H_{2}\right)  \leq v\left(  H_{2}\right)  c\left(  H_{2}\right)
\leq\left(  n-\left\vert M_{1}\right\vert -\left\vert H_{1}\right\vert
\right)  \left(  \frac{1}{4}+\frac{\alpha}{2}\right)  n.
\]
Thus, in view of
\[
e\left(  H_{2}\right)  =e\left(  G-M_{1}-H_{1}\right)  >\left(  \frac{1}%
{4}-\beta\right)  n^{2}-\frac{9n}{40}\left\vert M_{0}\right\vert
-n-\frac{\left\vert H_{1}\right\vert ^{2}}{2},
\]
and the previous inequality we have
\[
\left(  n-\left\vert M_{1}\right\vert -\left\vert H_{1}\right\vert \right)
\left(  \frac{1}{2}+\alpha\right)  n>\left(  \frac{1}{2}-2\beta\right)
n^{2}-\frac{9n}{20}\left\vert M_{0}\right\vert -2n-\left\vert H_{1}\right\vert
^{2},
\]
and so,%
\[
\left\vert H_{1}\right\vert ^{2}-\frac{1}{2}n\left\vert H_{1}\right\vert
+\left(  \alpha+2\beta\right)  n^{2}+2n\geq\alpha\left\vert H_{1}\right\vert
+\left(  \frac{1}{20}+\alpha\right)  \left\vert M_{0}\right\vert n+\left(
\frac{1}{2}+\alpha\right)  \left\vert K\right\vert >0.
\]
Solving the quadratic inequality with respect to $\left\vert H_{1}\right\vert
$ we see that%
\begin{equation}
\left\vert H_{1}\right\vert \geq\frac{1+\sqrt{1-16\left(  \alpha
+2\beta\right)  -32/n}}{4}n \label{lob}%
\end{equation}
or\
\[
\left\vert H_{1}\right\vert \leq\frac{1-\sqrt{1-16\left(  \alpha
+2\beta\right)  -32/n}}{4}n.
\]
Since
\begin{align*}
\frac{1-\sqrt{1-16\left(  \alpha+2\beta\right)  -32/n}}{4}  &  \leq
\frac{1-1+16\left(  \alpha+2\beta\right)  +32/n}{4}\\
&  =4\left(  \alpha+2\beta\right)  +\frac{32}{n}<1/3,
\end{align*}
we see that precisely (\ref{lob}) holds. From
\[
1>1-16\left(  \alpha+2\beta\right)  -\frac{32}{n}>0,
\]
we deduce that
\[
\sqrt{1-16\left(  \alpha+2\beta\right)  -32/n}>1-16\left(  \alpha
+2\beta\right)  -32/n,
\]
and so,
\begin{align*}
\left\vert H_{1}\right\vert  &  \geq\frac{1+\sqrt{1-16\left(  \alpha
+2\beta\right)  -32/n}}{4}n\geq\frac{1+1-16\left(  \alpha+2\beta\right)
-32/n}{4}n\\
&  =\left(  \frac{1}{2}-4\left(  \alpha+2\beta\right)  \right)  n-8\geq\left(
\frac{1}{2}-20\left(  \alpha+2\beta\right)  \right)  n.
\end{align*}
This, together with%
\[
\left\vert H_{2}\right\vert \leq n-\left\vert H_{1}\right\vert \leq\left(
\frac{1}{2}+20\left(  \alpha+2\beta\right)  \right)  n,
\]
completes the proof of (\ref{in0}).

To complete the proof of the theorem, we shall remove all low degree vertices
from $H_{1}\cup H_{2}.$ Letting
\[
M_{2}=\left\{  v:v\in V\left(  H_{1}\cup H_{2}\right)  ,\text{ }d_{H_{1}\cup
H_{2}}\left(  v\right)  \leq\frac{9}{20}n\right\}  ,
\]
we find that%
\begin{align*}
\left(  \frac{1}{2}-2\beta\right)  n^{2}  &  <2e\left(  G\right)  =\sum_{u\in
V\left(  G\right)  }d\left(  u\right)  =\sum_{u\in V\left(  H_{1}\cup
H_{2}\right)  \backslash M_{2}}d\left(  u\right)  +\sum_{u\in M_{2}}d\left(
u\right)  +\sum_{u\in M_{1}}d\left(  u\right) \\
&  <\left(  n-\left\vert M_{1}\right\vert -\left\vert M_{2}\right\vert
\right)  \left\vert H_{2}\right\vert +\frac{9}{20}\left\vert M_{2}\right\vert
n+\left\vert M_{1}\right\vert n\\
&  \leq\left(  n-\left\vert M_{1}\right\vert -\left\vert M_{2}\right\vert
\right)  \left(  \frac{1}{2}+20\left(  \alpha+2\beta\right)  \right)
n+\frac{9}{20}\left\vert M_{2}\right\vert n+\left\vert M_{1}\right\vert n\\
&  \leq\left(  \frac{1}{2}+20\left(  \alpha+2\beta\right)  \right)
n^{2}+\frac{1}{2}\left\vert M_{1}\right\vert n-\frac{1}{20}\left\vert
M_{2}\right\vert n\\
&  \leq\left(  \frac{1}{2}+20\left(  \alpha+2\beta\right)  \right)
n^{2}+20\left(  \alpha+\beta\right)  n^{2}-\frac{1}{20}\left\vert
M_{2}\right\vert n\\
&  =\left(  \frac{1}{2}+40\alpha+60\beta\right)  n^{2}-\frac{1}{20}\left\vert
M_{2}\right\vert n,
\end{align*}
and hence, $\left\vert M_{2}\right\vert \leq\left(  800\alpha+1240\beta
\right)  n.$ Setting
\[
M=M_{1}\cup M_{2}\text{, \ \ }G_{1}=H_{1}-M_{2},\text{ \ \ }G_{2}=H_{2}%
-M_{2}\text{,}%
\]
we see that
\begin{align*}
\left\vert M\right\vert  &  =\left\vert M_{1}\right\vert +\left\vert
M_{2}\right\vert \leq\left(  840\alpha+1280\beta\right)  n<840\left(
\alpha+2\beta\right)  n,\\
\left\vert G_{1}\right\vert  &  \geq\left\vert H_{1}\right\vert -\left\vert
M_{2}\right\vert \geq\left(  \frac{1}{2}-840\left(  \alpha+2\beta\right)
\right)  n,\\
\left\vert G_{2}\right\vert  &  \leq\left\vert H_{2}\right\vert \leq\left(
\frac{1}{2}+20\left(  \alpha+2\beta\right)  \right)  n,\\
\delta\left(  G_{1}\cup G_{2}\right)   &  \geq\frac{9}{20}n-\left\vert
M_{2}\right\vert \geq\left(  \frac{9}{20}-800\alpha-1240\beta\right)  n\\
&  \geq\left(  \frac{9}{20}-\frac{800}{10^{5}}-\frac{1240}{10^{5}}\right)
n>\frac{3}{7}n.
\end{align*}
Since
\[
\frac{3}{7}n>\frac{1}{2}\left(  \frac{1}{2}+20\left(  \alpha+2\beta\right)
\right)  n\geq\frac{1}{2}\max\left\{  \left\vert G_{2}\right\vert ,\left\vert
G_{1}\right\vert \right\}  ,
\]
it follows that $G_{1}$ and $G_{2}$ are connected, completing the proof.
\end{proof}

\subsection{Proof of Theorem \ref{cycth}}

Assume that $C_{l}\nsubseteq G$ for some $l\in\left[  3,\left\lceil \left(
1/2+\gamma\right)  n\right\rceil \right]  ;$ then Theorem \ref{ThBe} implies
\begin{equation}
c\left(  G\right)  \leq\left\lfloor \left(  1/2+\gamma\right)  n\right\rfloor
. \label{upc}%
\end{equation}

First we shall show that the assertion of the theorem holds for $n<\gamma
^{-1}/2.$ Indeed, by Theorem \ref{leEG}, for $n$ even, say $n=2k,$ we have
$\left\lfloor k+\gamma n\right\rfloor \geq c\left(  G\right)  >k,$
contradicting (\ref{upc}) for $n<\gamma^{-1}/2<\gamma^{-1}.$ Similarly, for
$n$ odd, say $n=2k+1,$ we have
\[
\left\lfloor k+\frac{1}{2}+\gamma n\right\rfloor \geq c\left(  G\right)  \geq
k+1,
\]
contradicting (\ref{upc}) for $n<\gamma^{-1}/2.$

In view of $n\geq\gamma^{-1}/2$ and (\ref{upc}), Theorem \ref{thDC}, with
$\alpha=\gamma$ and $\beta=0,$ implies that there exists $M\subset V\left(
G\right)  $ with $\left\vert M\right\vert <840\gamma n$ such that $G-M$
consists of two components $G_{1}$ and $G_{2}$ satisfying%
\begin{align*}
\left(  \frac{1}{2}-840\gamma\right)  n  &  <\left\vert G_{1}\right\vert
\leq\left\vert G_{2}\right\vert <\left(  \frac{1}{2}+20\gamma\right)  n\\
\delta\left(  G_{1}\right)   &  \geq3n/7,\text{ \ \ \ \ }\delta\left(
G_{1}\right)  \geq3n/7.
\end{align*}

From
\[
\frac{3n}{7}\geq\frac{1}{2}\left(  \left(  \frac{1}{2}+20\gamma\right)
n+1\right)  \geq\frac{1}{2}\left(  \left\vert G_{2}\right\vert +1\right)  ,
\]
and Lemma \ref{2hl}, we see that $G_{1}$ and $G_{2}$ are Hamiltonian connected.

Suppose there are two vertex disjoint paths $P\left(  u_{1},v_{1}\right)  $
and $P\left(  u_{2},v_{2}\right)  $ joining vertices from $G_{1}$ to vertices
from $G_{2},$ say
\begin{align*}
P\left(  v_{1},u_{1}\right)  \cap G_{1}  &  =\left\{  u_{1}\right\}  ,\text{
\ \ }P\left(  u_{1},v_{1}\right)  \cap G_{2}=\left\{  v_{1}\right\}  ,\\
P\left(  u_{2},v_{2}\right)  \cap G_{1}  &  =\left\{  u_{2}\right\}  ,\text{
\ \ }P\left(  u_{2},v_{2}\right)  \cap G_{2}=\left\{  v_{2}\right\}  .
\end{align*}
Let $Q_{1}\left(  u_{1},u_{2}\right)  $ and $Q_{2}\left(  v_{2},v_{1}\right)
$ be Hamiltonian paths within $G_{1}$ and $G_{2}.$ Then the length of the
cycle
\[
Q_{1}\left(  u_{1},u_{2}\right)  P\left(  u_{2},v_{2}\right)  Q_{2}\left(
v_{2},v_{1}\right)  P\left(  v_{1},u_{1}\right)
\]
is at least
\[
\left\vert G_{1}\right\vert +\left\vert G_{2}\right\vert =n-\left\vert
M\right\vert >\left(  1-840\gamma\right)  n>\left(  \frac{1}{2}+\gamma\right)
n,
\]
contradicting (\ref{upc}).

Therefore, no two vertex-disjoint paths join vertices from $G_{1}$ to vertices
from $G_{2}.$ By Menger's theorem, there exists a vertex $u\in V\left(
G\right)  $ separating $G_{1}$ and $G_{2}.$ Clearly, $V\left(  G_{1}\right)
\backslash\left\{  u\right\}  $ induces a connected subgraph in $G-u;$ let
$H_{1}$ be the component containing $V\left(  G_{1}\right)  \backslash\left\{
u\right\}  ,$ and $H_{2}$ be the union of the remaining components of $G-u$.
Observing that
\begin{align*}
\left(  \frac{1}{2}+840\gamma\right)  n  &  >n-\left\vert G_{2}\right\vert
\geq\left\vert H_{1}\right\vert \geq\left\vert G_{1}\right\vert -1>\left(
\frac{1}{2}-900\gamma\right)  n,\\
\left(  \frac{1}{2}+840\gamma\right)  n  &  \geq n-\left\vert G_{1}\right\vert
\geq\left\vert H_{2}\right\vert \geq\left\vert G_{2}\right\vert -1>\left(
\frac{1}{2}-900\gamma\right)  n,
\end{align*}
we complete the proof.\hfill$\square$

\section{\label{Ram}Ramsey type results}

Theorem \ref{th3par}, presented in the beginning of this section, is
essentially a stability result of Tur\'{a}n type. However, it is placed in
this section, since it is the main tool to derive Theorem \ref{arrth} - a
distinctive Ramsey type result.

\begin{theorem}
\label{th3par} Let $0<\alpha<5\times10^{-6},$ $0\leq\beta\leq\alpha/25,$ and
$n\geq\alpha^{-1}.$ If $G=G\left(  n\right)  $ is a graph with $e\left(
G\right)  >\left(  1/4-\beta\right)  n^{2}$, then one of the following
conditions hold:

(i) $C_{t}\subset G$ for every $t\in\left[  3,\left\lceil \left(
1/2+\alpha\right)  n\right\rceil \right]  $;

(ii) there exists a partition $V\left(  G\right)  =V_{0}\cup V_{1}\cup V_{2}$
such that%
\begin{align}
\left\vert V_{0}\right\vert  &  <2000\alpha n,\label{i1}\\
\left(  \frac{1}{2}-10\sqrt{\alpha+\beta}\right)  n  &  <\left\vert
V_{1}\right\vert \leq\left\vert V_{2}\right\vert <\left(  \frac{1}{2}%
+10\sqrt{\alpha+\beta}\right)  n,\label{i2}\\
\delta\left(  G-V_{0}\right)   &  \geq2n/5, \label{i3}%
\end{align}
and either
\[
E\left(  G-V_{0}\right)  \subset V_{1}^{\left(  2\right)  }\cup V_{2}^{\left(
2\right)  }\text{ \ \ \ or \ \ \ }E\left(  G-V_{0}\right)  \subset V_{1}\times
V_{2}.
\]

\end{theorem}

\begin{proof}
Setting
\[
M=\left\{  v:v\in V\left(  G\right)  ,\text{ }d\left(  v\right)  \leq\frac
{9n}{20}\right\}  ,
\]
our first goal is to prove that
\begin{equation}
\left\vert M\right\vert <20\left(  \alpha+2\beta\right)  n. \label{ubm1}%
\end{equation}
Indeed, assume for contradiction that $\left\vert M\right\vert \geq20\left(
\alpha+2\beta\right)  n>24\beta n$ and select $M_{0}\subset\left\vert
M\right\vert $ with $\left\vert M_{0}\right\vert =\left\lceil 24\beta
n\right\rceil .$ We shall show that
\begin{equation}
e\left(  G-M_{0}\right)  >\frac{1}{4}\left(  n-\left\vert M_{0}\right\vert
\right)  ^{2}. \label{emin}%
\end{equation}
Indeed, otherwise we have
\begin{align*}
\left(  \frac{1}{4}-\beta\right)  n^{2}  &  <e\left(  G\right)  \leq e\left(
G-M_{0}\right)  +\sum_{u\in M_{0}}d\left(  u\right)  \leq e\left(
G-M_{0}\right)  +\frac{9n}{20}\left\vert M_{0}\right\vert \\
&  \leq\frac{1}{4}\left(  n-\left\vert M_{0}\right\vert \right)  ^{2}%
+\frac{9n}{20}\left\vert M_{0}\right\vert =\frac{1}{4}n^{2}-\frac{1}%
{20}n\left\vert M_{0}\right\vert +\frac{\left\vert M_{0}\right\vert ^{2}}{4},
\end{align*}
and so,
\[
20\beta n^{2}-n\left\vert M_{0}\right\vert +5\left\vert M_{0}\right\vert
^{2}\geq0.
\]
Solving this quadratic inequality with respect to $\left\vert M_{0}\right\vert
,$ we see that either%
\[
\left\vert M_{0}\right\vert \leq\frac{1-\sqrt{1-400\beta}}{10}n<24\beta
n\leq\left\lceil 24\beta n\right\rceil
\]
or
\[
\left\vert M_{0}\right\vert \geq\frac{1+\sqrt{1-400\beta}}{10}n>\frac{1}%
{10}n>24\beta n+1>\left\lceil 24\beta n\right\rceil .
\]
Since both inequalities contradict our choice of $M_{0},$ inequality
(\ref{emin}) holds.

Note that, in view of
\[
\left(  \alpha-12\beta-48\alpha\beta\right)  n\geq\left(  \alpha-\frac{12}%
{25}\alpha-\frac{2\cdot48\alpha}{25\cdot100,000}\right)  n\geq\frac{51}%
{100}\alpha n>\frac{1}{2}+2\alpha,
\]
we have
\begin{align*}
\left(  \frac{1}{2}+2\alpha\right)  \left\vert G-M_{0}\right\vert  &  =\left(
\frac{1}{2}+2\alpha\right)  \left(  n-\left\lceil 24\beta n\right\rceil
\right)  \geq\left(  \frac{1}{2}+2\alpha\right)  \left(  n-24\beta n-1\right)
\\
&  \geq\left(  \frac{1}{2}+\alpha\right)  n+\left(  \alpha-12\beta
-48\alpha\beta\right)  n-\left(  \frac{1}{2}+2\alpha\right)  \geq\left(
\frac{1}{2}+\alpha\right)  n.
\end{align*}
Hence, if $C_{t}\subset G-M_{0}$ for every $t\in\left[  3,\left\lceil \left(
1/2+2\alpha\right)  \left\vert G-M_{0}\right\vert \right\rceil \right]  ,$ we
see that \emph{(i)} holds. Thus, $C_{t}\nsubseteq G-M_{0}$ for some
$t\in\left[  3,\left\lceil \left(  1/2+2\alpha\right)  \left\vert
G-M_{0}\right\vert \right\rceil \right]  .$ Applying Theorem \ref{thDC} to the
graph $G-M_{0}$ with $\alpha^{\prime}=2\alpha,$ $\beta^{\prime}=0,$ it follows
that there exists a $M_{1}\subset V\left(  G-M_{0}\right)  $ such that
$G-M_{0}-M_{1}=G_{1}\cup G_{2},$ where $G_{1}$ and $G_{2}$ are vertex-disjoint
graphs satisfying
\begin{align*}
\left\vert M_{1}\right\vert  &  <1680\alpha\left(  n-\left\vert M_{0}%
\right\vert \right)  <1800\alpha\left(  n-\left\vert M_{0}\right\vert \right)
,\\
\left(  \frac{1}{2}-1800\alpha\right)  \left(  n-\left\vert M_{0}\right\vert
\right)   &  <\left\vert G_{1}\right\vert \leq\left\vert G_{2}\right\vert
<\left(  \frac{1}{2}+40\alpha\right)  \left(  n-\left\vert M_{0}\right\vert
\right)  ,\\
\delta\left(  G_{1}\right)   &  \geq\frac{3}{7}\left(  n-\left\vert
M_{0}\right\vert \right)  ,\text{ \ \ \ \ }\delta\left(  G_{2}\right)
\geq\frac{3}{7}\left(  n-\left\vert M_{0}\right\vert \right)  .
\end{align*}
Setting
\[
V_{0}=M_{0}\cup M_{1},\text{ \ \ }V_{1}=V\left(  G_{1}\right)  ,\text{
\ \ }V_{2}=V\left(  G_{2}\right)  ,
\]
we first note that $E\left(  G-V_{0}\right)  \subset V_{1}^{\left(  2\right)
}\cup V_{2}^{\left(  2\right)  }.$ We shall prove that this selection of
$V_{0},$ $V_{1},$ and $V_{2}$ satisfies \emph{(ii).} To this end we have to
derive inequalities (\ref{i1}), (\ref{i2}) and (\ref{i3}). Inequality
(\ref{i1}) follows from
\[
\left\vert V_{0}\right\vert \leq\left\vert M_{0}\right\vert +1800\alpha\left(
n-\left\vert M_{0}\right\vert \right)  \leq\left\lceil 24\beta n\right\rceil
+1800\alpha n<24\beta n+1+1800\alpha n<2000\alpha n.
\]
Our next goal is to prove (\ref{i2}). Note that
\begin{align*}
\left\vert G_{1}\right\vert  &  \geq\left(  \frac{1}{2}-1800\alpha\right)
\left(  n-\left\vert M_{0}\right\vert \right)  >\left(  \frac{1}{2}%
-1800\alpha\right)  \left(  n-24\beta n-1\right) \\
&  \geq\frac{n-1}{2}-12\beta n-1800\alpha n\geq\left(  200\alpha
-12\beta\right)  n-\frac{1}{2}+\left(  \frac{1}{2}-2000\alpha\right)  n\\
&  \geq\left(  \frac{1}{2}-2000\alpha\right)  n\geq\left(  \frac{1}{2}%
-10\sqrt{\alpha}\right)  n\geq\left(  \frac{1}{2}-10\sqrt{\alpha+\beta
}\right)  n.
\end{align*}
Since
\[
\left\vert G_{2}\right\vert <n-\left\vert G_{1}\right\vert \leq\left(
\frac{1}{2}+10\sqrt{\alpha+\beta}\right)  n,
\]
inequality (\ref{i2}) follows. Finally, (\ref{i3}) follows from
\[
\delta\left(  G-V_{0}\right)  \geq\frac{3}{7}\left(  n-\left\vert
M_{0}\right\vert \right)  =\frac{3}{7}\left(  n-\left\lceil 24\beta
n\right\rceil \right)  \geq\frac{3}{7}\left(  n-24\beta n-1\right)  >\frac
{2}{5}n.
\]

This completes the proof of the theorem if (\ref{ubm1}) fails. Thus,
hereafter, we shall assume that (\ref{ubm1}) holds. Set $V_{0}=M,$
$G_{0}=G-V_{0},$ and observe that
\begin{align}
\left\vert V_{0}\right\vert  &  =\left\vert M\right\vert \leq20\left(
\alpha+2\beta\right)  n<2000\alpha n,\label{in1.1}\\
\delta\left(  G_{0}\right)   &  \geq\frac{9n}{20}-\left\vert M\right\vert
=\left(  \frac{9}{20}-20\left(  \alpha+2\beta\right)  \right)  n>\frac{2}{5}n.
\label{in1.3}%
\end{align}

\emph{Case 1: }$G_{0}$\emph{ is bipartite}

Write $V_{1}$ and $V_{2}$ for the vertex classes of $G_{0},$ and let say,
$\left\vert V_{1}\right\vert \leq\left\vert V_{2}\right\vert $. We see that
$E\left(  G_{0}\right)  \subset V_{1}\times V_{2},$ also (\ref{i1}) and
(\ref{i3})\ hold in view of (\ref{in1.1}) and (\ref{in1.3}), so to finish the
proof, we need to prove inequalities (\ref{i2}). Since
\begin{align*}
e\left(  G_{0}\right)   &  \geq\left(  \frac{1}{4}-\beta\right)
n^{2}-\left\vert V_{0}\right\vert n\geq\left(  \frac{1}{4}-20\alpha
-41\beta\right)  n^{2}\\
&  >\left(  \frac{1}{4}-\left(  10\sqrt{\alpha+\beta}\right)  ^{2}\right)
n^{2},
\end{align*}
selecting $x$ so that%
\[
\left\vert V_{1}\right\vert =\left(  \frac{1}{2}-x\right)  \left(  \left\vert
V_{1}\right\vert +\left\vert V_{2}\right\vert \right)  ,\text{ \ \ }\left\vert
V_{2}\right\vert =\left(  \frac{1}{2}+x\right)  \left(  \left\vert
V_{1}\right\vert +\left\vert V_{2}\right\vert \right)  ,
\]
we deduce that%
\[
\left(  \frac{1}{4}-x^{2}\right)  \geq\left(  \frac{1}{4}-\left(
10\sqrt{\alpha+\beta}\right)  ^{2}\right)  ,
\]
and,%
\[
\left\vert V_{2}\right\vert =\left(  \frac{1}{2}+x\right)  \left(  \left\vert
V_{1}\right\vert +\left\vert V_{2}\right\vert \right)  <\left(  \frac{1}%
{2}+10\sqrt{\alpha+\beta}\right)  n.
\]
This inequality implies in turn
\[
\left\vert V_{1}\right\vert \geq n-\left\vert V_{2}\right\vert >\left(
\frac{1}{2}-10\sqrt{\alpha+\beta}\right)  n,
\]
completing the proof in this case.

\emph{Case 2: }$\kappa\left(  G_{0}\right)  \leq1$

Let $K$ be a cutset of $G_{0}$ with $\left\vert K\right\vert \leq1.$ Since
$\delta\left(  G_{0}-K\right)  >2n/5-1>n/3,$ the graph $G_{0}-K$ has exactly
two components - $G_{1}$ and $G_{2}.$ Let $V_{0}=M\cup K,$ $V_{1}=V\left(
G_{1}\right)  ,$ $V_{2}=V\left(  G_{2}\right)  ;$ assume $\left\vert
V_{1}\right\vert \leq\left\vert V_{2}\right\vert $ and observe that $E\left(
G_{0}\right)  \subset V_{1}^{\left(  2\right)  }\cup V_{2}^{\left(  2\right)
}.$ Clearly $\left\vert V_{0}\right\vert \leq20\left(  \alpha+2\beta\right)
n+1\leq2000\alpha n,$ so (\ref{i1}) holds. From
\[
\delta\left(  G-V_{0}\right)  >\frac{9n}{20}-\left\vert M\right\vert
-1>\frac{3}{8}n>\frac{n-\left\vert V_{1}\right\vert }{2}>\frac{1}{2}\left\vert
V_{2}\right\vert ,
\]
we see first, that(\ref{i3}) holds, and second, that $G_{2}$ is Hamiltonian.
From Theorem \ref{thBFG} it follows that $C_{t}\subset G_{2}$ for every
$t\in\left[  3,\left\vert V_{2}\right\vert \right]  .$ This completes the
proof of the theorem if
\begin{equation}
\left\vert V_{2}\right\vert \geq\left(  \frac{1}{2}+5\sqrt{\alpha+2\beta
}\right)  n, \label{in4}%
\end{equation}
since then $\left\vert V_{2}\right\vert >\left(  1/2+\alpha\right)  n,$ and so
\emph{(i)} holds.

Assume that (\ref{in4}) fails. Then
\[
\left\vert V_{2}\right\vert <\left(  \frac{1}{2}+5\sqrt{\alpha+2\beta}\right)
n<\left(  \frac{1}{2}+10\sqrt{\alpha+\beta}\right)  ,
\]
and so
\[
\left\vert V_{1}\right\vert >n-\left\vert V_{2}\right\vert >\left(  \frac
{1}{2}-10\sqrt{\alpha+\beta}\right)  n.
\]
Thus (\ref{i2}) holds, completing the proof of \emph{(ii) }in this case.

\emph{Case 3: }$G_{0}$\emph{ is }$2$\emph{-connected and nonbipartite}

In this case we shall show that \emph{(i)} holds.\emph{ }Since $\delta\left(
G_{0}\right)  >2n/5,$ Dirac's theorem implies that $c\left(  G\right)
\geq2\delta\left(  G_{0}\right)  >4n/5>\left\lceil \left(  1/2+\alpha\right)
n\right\rceil .$ Now, Theorem \ref{thBFG} implies that $C_{t}\subset G_{0}$
for all $t\in\left[  3,\left\lceil \left(  1/2+\alpha\right)  n\right\rceil
\right]  ,$ completing the proof.
\end{proof}

\subsection{Proof of Theorem \ref{cycth1}}

We precede the proof of Theorem \ref{cycth1} by a simple lemma whose idea goes
back to \cite{FaSc74}. The present version of the lemma emerged from recent
conversations with Ingo Schiermeyer, Linda Lesniak, and Ralph Faudree, to whom
we are grateful. The lemma helped enhance considerably an earlier version of
Theorem \ref{cycth1}.

\begin{lemma}
\label{Le4}Let $G$ be a Hamiltonian graph of order $2n$ such that
$C_{2n-1}\nsubseteq G$ and $C_{2n-1}\nsubseteq\overline{G}.$ Then there exists
a partition $V\left(  G\right)  =U_{1}\cup U_{2}$ such that $\left\vert
U_{1}\right\vert =\left\vert U_{2}\right\vert =n$ and $U_{1},$ $U_{2}$ are
independent. Moreover, there exists a vertex $u\in V\left(  G\right)  $ such
that $G-u=K_{n,n-1}$.
\end{lemma}

\begin{proof}
Assume $v_{1},v_{2},...,v_{2n}$ are the vertices of $G$ listed along the
Hamiltonian cycle of $G.$ Observe that $\left(  v_{1},v_{3},...,v_{2n-1}%
,v_{1}\right)  $ and $\left(  v_{2},v_{4},...,v_{2n},v_{2}\right)  $ are
cycles of length $n$ in $\overline{G}.$ Our first goal is to show that the
sets $U_{1}=\left\{  v_{1},v_{3},...,v_{2n-1}\right\}  $ and $U_{2}=\left\{
v_{2},v_{4},...,v_{2n}\right\}  $ are independent. Assume for contradiction
that this is not true and let say $v_{1}v_{2k+1}\in E\left(  G\right)  .$ Then
$v_{3}v_{2k+2}\notin E\left(  G\right)  $ since otherwise,
\[
\left(  v_{3},v_{4},...,v_{2k+1},v_{1},v_{2n},v_{2n-1}...,v_{2k+2}\right)
\]
is a cycle of length $2n-1$ in $G.$ Likewise, $v_{2n-1}v_{2k}\notin E\left(
G\right)  .$ Then
\[
\left(  v_{3},v_{5},...,v_{2n-1},v_{2k},v_{2k-2},...,v_{2k+2}\right)
\]
is a cycle of length $2n-1$ in $\overline{G},$ a contradiction. Therefore
$\overline{G}\left[  U_{1}\right]  $ and $\overline{G}\left[  U_{2}\right]  $
are complete graphs. Since $C_{2n-1}\nsubseteq\overline{G},$ we see that
$E_{\overline{G}}\left(  U_{1},U_{2}\right)  $ contains no disjoint edges and
therefore is a (possibly empty) star. Taking $u$ to be the center of this star
we complete the proof.
\end{proof}

\begin{proof}
[\textbf{Proof of Theorem \ref{cycth1}}]Recall that if $k$ is an integer, then
for every $2$-coloring $E\left(  K_{3k-1}\right)  =E\left(  R\right)  \cup
E\left(  B\right)  ,$ either $C_{2k}\subset R$ or $C_{2k}\subset B$ (e.g., see
\cite{FaSc74}). Since for $\beta\leq1/2,$ $\left\lfloor \left(  2-\beta
\right)  2k\right\rfloor \geq3k>3k-1,$ we see that the assertion holds
immediately for even $n$. Let $n$ be odd, say $n=2k+1,$ set $V=V\left(
K_{\left\lfloor \left(  2-\beta\right)  \left(  2k+1\right)  \right\rfloor
}\right)  ,$ and assume $E\left(  K_{\left\lfloor \left(  2-\beta\right)
\left(  2k+1\right)  \right\rfloor }\right)  =E\left(  R\right)  \cup E\left(
B\right)  $ is a $2$-coloring with $C_{2k+1}\nsubseteq R$ and $C_{2k+1}%
\nsubseteq B$. From%
\[
\left(  2-\beta\right)  \left(  2k+1\right)  \geq\left(  2-\frac{k}%
{2k+1}\right)  \left(  2k+1\right)  =3k+2
\]
we see that, up to color, $C_{2k+2}\subset B$. By the assumption of the
theorem and Lemma \ref{Le4} it follows that there exist $W_{1},W_{2}\subset V$
such that $\left\vert W_{1}\right\vert =k,$ $\left\vert W_{2}\right\vert
=k+1,$ $W_{1}$ and $W_{2}$ induce complete graphs in $R,$ and $K_{k,k+1}$ in
$B.$ Note that for all $u\in V\backslash\left(  W_{1}\cup W_{2}\right)  $
either $\Gamma_{B}\cap W_{1}=\varnothing$ or $\Gamma_{B}\cap W_{2}%
=\varnothing,$ as otherwise $C_{2k+1}\subset B.$ Set
\begin{align*}
X_{1}  &  =\left\{  u:u\in V\backslash\left(  W_{1}\cup W_{2}\right)  \text{
and }\Gamma_{B}\cap W_{1}=\varnothing\right\}  ,\\
X_{2}  &  =V\backslash\left(  W_{1}\cup W_{2}\cup U_{1}\right)  ,\\
V_{1}  &  =X_{1}\cup W_{1},\text{ \ \ }V_{2}=X_{2}\cup W_{2},
\end{align*}
and note that $X_{1}\times W_{1}\subset R$ and $X_{2}\times W_{2}\subset R.$
At this stage it is not difficult to check immediately that the assertion of
the theorem holds for $k=2,$ so in the sequel we shall assume that $k\geq3.$

If there exist two disjoint edges $v_{1}u_{1},u_{2}v_{2}\in E_{R}\left(
V_{1},V_{2}\right)  ,$ then $C_{t}\subset R$ for any odd $t\in\left[
7,2k+1\right]  .$ Hence, $E_{R}\left(  V_{1},V_{2}\right)  $ is a (possibly
empty) star; let $u$ be its center or any other vertex if $E_{R}\left(
V_{1},V_{2}\right)  $ is empty; set $U_{1}=V_{1}\backslash\left\{  u\right\}
,$ $U_{2}=V_{2}\backslash\left\{  u\right\}  $. Then $U_{1}\times U_{2}\subset
E\left(  B\right)  $ and hence, $E_{B}\left(  U_{1}\right)  =E_{B}\left(
U_{2}\right)  =\varnothing,$ as otherwise $C_{2k+1}\subset B.$ To prove
inequalities (\ref{thin}), we shall assume that $\left\vert U_{1}\right\vert
\leq\left\vert U_{2}\right\vert .$ This implies that $\left\vert
U_{2}\right\vert \leq2k,$ and so
\[
\left\vert U_{1}\right\vert =\left\lfloor \left(  2-\beta\right)  \left(
2k+1\right)  \right\rfloor -2k-1>\left(  2-\beta\right)  \left(  2k+1\right)
-2k-2=\left(  1-\beta\right)  \left(  2k+1\right)  -1,
\]
completing the proof.
\end{proof}

\subsection{Proof of Theorem \ref{arrth}}

For convenience we shall rephrase the Theorem \ref{arrth} in terms of
$3$-colorings of $K_{2n-1}.$

\begin{theorem}
Let the edges of $K_{2n-1}$ be $3$-colored, i.e., $E\left(  K_{2n-1}\right)  $
be partitioned as $E\left(  K_{2n-1}\right)  =E\left(  R\right)  \cup E\left(
B\right)  \cup E\left(  Y\right)  ,$ where $R,$ $B,$ and $Y$ are graphs with
$V\left(  R\right)  =V\left(  B\right)  =V\left(  Y\right)  =\left[
2n-1\right]  .$ Let the minimum degree $\delta\left(  R\cup B\right)  $
satisfies $\delta\left(  R\cup B\right)  >\left(  2-10^{-6}\right)  n.$ Then,
if $n$ is sufficiently large, either $C_{t}\subset R$ for all $t\in\left[
3,n\right]  $ or $C_{t}\subset B$ for all $t\in\left[  3,n\right]  .$
\end{theorem}

\begin{proof}
Set for brevity
\begin{align}
c  &  =10^{-6},\nonumber\\
\beta &  =c/8=10^{-6}/8,\label{defb}\\
\alpha &  =25\beta=10^{-4}/32, \label{defa}%
\end{align}
and assume, without loss of generality, that $e\left(  R\right)  \geq e\left(
B\right)  .$ Hence, from%
\[
e\left(  R\right)  +e\left(  B\right)  \geq\left(  2-10^{-6}\right)  \left(
2n-1\right)  n\geq\left(  \frac{1}{2}-2\beta\right)  \left(  2n-1\right)
^{2},
\]
we see that%
\[
e\left(  R\right)  \geq\left(  \frac{1}{4}-\beta\right)  \left(  2n-1\right)
^{2}.
\]
According to Theorem \ref{th3par}, one of the following conditions hold:

\emph{(i)} $C_{t}\subset R$ for every $t\in\left[  3,\left\lceil \left(
1/2+\alpha\right)  \left(  2n-1\right)  \right\rceil \right]  $;

\emph{(ii)} there exists a partition $\left[  2n-1\right]  =V_{0}\cup
V_{1}\cup V_{2}$ such that
\begin{align*}
\left\vert V_{0}\right\vert  &  <2000\alpha\left(  2n-1\right)  ,\\
\left(  \frac{1}{2}-10\sqrt{\alpha+\beta}\right)  \left(  2n-1\right)   &
<\left\vert V_{1}\right\vert \leq\left\vert V_{2}\right\vert <\left(  \frac
{1}{2}+10\sqrt{\alpha+\beta}\right)  \left(  2n-1\right)  ,
\end{align*}
and either
\[
E\left(  R-V_{0}\right)  \subset V_{1}^{\left(  2\right)  }\cup V_{2}^{\left(
2\right)  }\text{ \ \ \ or \ \ \ }E\left(  R-V_{0}\right)  \subset V_{1}\times
V_{2}.
\]

If \emph{(i)}\ holds, there is nothing to prove, so we shall assume that
\emph{(ii)} holds. Then, in view of (\ref{defa}), (\ref{defb}), and%
\begin{align*}
2000\alpha &  =\frac{2000\cdot25}{8\cdot10^{6}}=\frac{1}{16\cdot10^{2}}%
=\frac{1}{160},\\
10\sqrt{\alpha+\beta}  &  =10\sqrt{\frac{26}{8\cdot10^{6}}}<\frac{1}{50},
\end{align*}
we find that
\begin{align*}
\left\vert V_{0}\right\vert  &  <\frac{1}{160}\left(  2n-1\right)  ,\\
\frac{12}{25}\left(  2n-1\right)   &  <\left\vert V_{1}\right\vert
\leq\left\vert V_{2}\right\vert <\frac{13}{25}\left(  2n-1\right)  .
\end{align*}

Assume $E\left(  R-V_{0}\right)  \subset V_{1}^{\left(  2\right)  }\cup
V_{2}^{\left(  2\right)  }$. We shall prove that, then $E\left(
B-V_{0}\right)  \subset V_{1}\times V_{2}.$ We clearly have
\begin{align*}
\delta\left(  B-V_{0}\right)   &  \geq\left\vert V_{1}\right\vert
-\Delta\left(  Y\right)  \geq\frac{12}{25}\left(  2n-1\right)  -\left(
\left(  2n-2\right)  -\delta\left(  R\cup B\right)  \right) \\
&  >\frac{12}{25}\left(  2n-1\right)  +2-cn>\left(  \frac{12}{25}-c\right)
\left(  2n-1\right)  \geq\frac{1}{2}\left(  \frac{13}{25}\left(  2n-1\right)
\right)  +1.
\end{align*}
Lemma \ref{ecl} implies that $C_{t}\subset B-V_{0}$ for all even $t\in\left[
4,2\left(  2\delta\left(  B-V_{0}\right)  -\left\vert V_{1}\right\vert
-1\right)  \right]  .$ Moreover, if $E\left(  B\left(  V_{1}\right)  \right)
\cup E\left(  B\left(  V_{2}\right)  \right)  \neq\varnothing,$ then obviously
$C_{t}\subset B-V_{0}$ for all odd $t\in\left[  3,2\left(  2\delta\left(
B-V_{0}\right)  -\left\vert V_{1}\right\vert \right)  \right]  .$ Since
\begin{align}
2\left(  2\delta\left(  B-V_{0}\right)  -\left\vert V_{1}\right\vert
-1\right)   &  \geq2\left(  2\left(  \frac{12}{25}-c\right)  \left(
2n-1\right)  -\frac{12}{25}\left(  2n-1\right)  -1\right) \label{ubn1}\\
&  \geq\left(  \left(  \frac{24}{25}-4c\right)  \left(  2n-1\right)
-1\right)  \geq n,\nonumber
\end{align}
the proof is completed. Hence, $E\left(  B\left(  V_{1}\right)  \right)  \cup
E\left(  B\left(  V_{2}\right)  \right)  =\varnothing,$ implying that
$E\left(  B-V_{0}\right)  \subset V_{1}\times V_{2}.$

Now, suppose that there exists $u\in V_{0}$ such that $\Gamma_{B}\left(
u\right)  \cap V_{1}\neq\varnothing$ and $\Gamma_{B}\left(  u\right)  \cap
V_{2}\neq\varnothing.$ Select $x\in\Gamma_{B}\left(  u\right)  \cap V_{1},$
$y\in\Gamma_{B}\left(  u\right)  \cap V_{2}$ and note that Lemma \ref{ecl}
implies that $B-V_{0}$ contains an $xy$-path of length $t$ for every odd
$t\in\left[  3,2\left(  2\delta\left(  B-V_{0}\right)  -\left\vert
V_{1}\right\vert -1\right)  \right]  .$ In view of (\ref{ubn1}), $C_{t}\subset
B-V_{0}$ for all odd $t\in\left[  3,n\right]  $, completing the proof.
Therefore, for every $u\in V_{0},$ either $\Gamma_{B}\left(  u\right)  \cap
V_{1}=\varnothing,$ or $\Gamma_{B}\left(  u\right)  \cap V_{2}=\varnothing.$
Set
\[
W_{1}=\left\{  u:u\in V_{0},\Gamma_{B}\left(  u\right)  \cap V_{1}%
=\varnothing\right\}  ,\text{ \ \ }W_{2}=V_{0}\backslash W_{1}.
\]
and let say $\left\vert W_{1}\cup V_{1}\right\vert \geq\left\vert W_{2}\cup
V_{2}\right\vert .$ This implies $\left\vert W_{1}\cup V_{1}\right\vert \geq
n.$ Let $t\in\left[  3,n\right]  .$ If $t\leq\left\vert V_{1}\right\vert ,$
then $C_{t}\subset R\left(  V_{1}\right)  \subset R.$ If $t>\left\vert
V_{1}\right\vert ,$ select $W_{1}^{\prime}\subset W_{1}$ so that $\left\vert
W_{1}^{\prime}\cup V_{1}\right\vert =t.$ Note that $V_{1}\times W_{1}^{\prime
}\subset E\left(  R\right)  \cup E\left(  Y\right)  $ and so for every $u\in
W_{1}^{\prime}\cup V_{1}$ we have
\begin{align*}
\left\vert \Gamma_{R}\left(  u\right)  \cap\left(  W_{1}^{\prime}\cup
V_{1}\right)  \right\vert  &  \geq\min\left\{  \left(  \left\vert
V_{1}\right\vert -\Delta\left(  Y\right)  \right)  ,\left\vert V_{1}%
\right\vert -1\right\} \\
&  \geq\frac{12}{25}\left(  2n-1\right)  -10^{-6}n-1\geq\frac{1}{2}n>\frac
{1}{2}t.
\end{align*}
Therefore, $R\left(  W_{1}^{\prime}\cup V_{1}\right)  $ is Hamiltonian, i.e.,
$C_{t}\subset R,$ completing the proof.
\end{proof}

.

\end{document}